\documentclass[11pt]{article}
\usepackage{extpfeil}
\usepackage[all]{xy}
\usepackage{amssymb} 
\usepackage[margin=1in]{geometry}
\usepackage{euscript}
\usepackage{mathtools}
\usepackage{amsthm}
\usepackage{url}

\numberwithin{equation}{section}
\theoremstyle{plain}
 \newtheorem{thm}[equation]{Theorem}
 \newtheorem{prop}[equation]{Proposition}
 \newtheorem{lem}[equation]{Lemma}
 \newtheorem{cor}[equation]{Corollary}

\theoremstyle{definition}
 \newtheorem{defn}[equation]{Definition}
 \newtheorem{exm}[equation]{Example}
\theoremstyle{remark}
 \newtheorem{rem}[equation]{Remark}

\DeclareMathAlphabet{\mathbbmsl}{U}{bbm}{m}{sl}

\def\on{\operatorname}
\renewcommand{\le}{\leqslant}
\renewcommand{\ge}{\geqslant}
\renewcommand{\setminus}{\smallsetminus}
\newcommand{\RR}{\mathbb{R}}
\newcommand{\m}{\mathfrak{m}}
\newcommand{\fC}{\EuScript C}
\newcommand{\ZZ}{\mathbb{Z}}
\def\AA{\mathbb{A}}
\newcommand{\Z}{\mathbb{Z}}

\newcommand{\R}{\mathbb{R}}
\newcommand{\C}{\mathbb{C}}
\newcommand{\CC}{\mathbb{C}}
\newcommand{\cC}{ {\EuScript C}}
\newcommand{\PP}{\mathbb{P}}
\newcommand \id {{\operatorname{id}}}
\newcommand\MF{{\operatorname{MF}}}
\newcommand\Hom{{\operatorname{Hom}}}

\newcommand\Ob{{\operatorname{Ob}}}

\newcommand\Spec{{\operatorname{Spec}}}

\newcommand{\into}{\hookrightarrow}

\newcommand{\K}{K}
\newcommand{\rK}{\widetilde{K}}
\newcommand{\GL}{\operatorname{GL}}
\newcommand{\BGL}{\operatorname{BGL}}
\newcommand{\BU}{\operatorname{BU}}
\newcommand{\EGL}{\operatorname{EGL}}
\newcommand{\charr}{\operatorname{char}}

\newcommand{\coker}{\operatorname{coker}}

\newcommand{\gr}{\operatorname{gr}}

\newcommand{\alg}{\operatorname{alg}}
\newcommand{\qgr}{\operatorname{qgr}}

\newcommand{\Db}{\operatorname{D}^b}

\newcommand{\topp}{\operatorname{top}}

\newcommand{\Proj}{\operatorname{Proj}}

\newcommand{\coh}{\operatorname{coh}}

\newcommand{\Sp}{\operatorname{Sp}}
\newcommand{\OO}{\mathcal{O}}
\newcommand{\RHom}{\operatorname{\bf{R}\!\operatorname{\underline{Hom}}}}

\newcommand{\sX}{\EuScript X}
\newcommand{\sU}{\EuScript U}

\newcommand{\sZ}{\EuScript Z}

\def\darrow#1#2{\xtofrom[#2]{#1}}
\newcommand{\sg}{\operatorname{sg}}

\newcommand{\Aff}{\on{Aff}_\C}
\newcommand{\Top}{\EuScript Top}
\newcommand{\M}{\EuScript M}

\def\lra{\longrightarrow}
\def\llra{\longleftrightarrow}
\def\dg{\on{dg}}
\def\Cat{\EuScript Cat}
\def\Catdg{{\Cat}_{\dg}}
\def\Lqe{\on{L}_{\on{qe}}(\Catdg(\k))}
\def\Lqec{\on{L}_{\on{qe}}(\Catdg(\C))}

\def\gr{\on{gr}\!}
\def\grn{\on{gr}_{\ge 0}\!}
\def\qgr{\on{qgr}\!}
\def\tor{\on{tor}\!}
\def\Perf{\on{\EuScript Perf}\!}
\def\k{\mathbf{k}}
\def\D{\EuScript D}
\def\Ch{\EuScript Ch}
\def\A{\EuScript A}

\def\Ds{\D^{\on{sg}}}
\def\L{\on{\bf L}\!}
\def\R{\on{\bf R}\!}

\title{Topological $K$-theory of equivariant singularity categories}

\author{Michael K. Brown\footnote{Hausdorff Center for Mathematics, Villa Maria, Endenicher Allee 62, 53115 Bonn, Germany, email: mbrown@math.uni-bonn.de, dyckerho@math.uni-bonn.de.}  \text{ }and Tobias Dyckerhoff\footnotemark[\value{footnote}]}

\begin{document}
\vspace{18mm} \setcounter{page}{1} \thispagestyle{empty}

\maketitle
\begin{abstract}
We study the topological $K$-theory spectrum of the dg singularity category associated to a weighted projective complete intersection. We calculate the topological $K$-theory of the dg singularity category of a weighted projective hypersurface in terms of its affine Milnor fiber and monodromy operator, and, as an application, we obtain a lift of the Atiyah-Bott-Shapiro construction to the level of spectra.
\end{abstract}
\tableofcontents

\section{Introduction}
Let $F$ be a field. The seminal work \cite{segal:categories} provides a beautiful method for
constructing the algebraic $K$-theory spectrum of $F$ from the sum operation on the category $\cC$
of finite-dimensional $F$-vector spaces. Namely, the sum equips the classifying space
\[
	\cC_1 := \coprod_{n \ge 0} \BGL(n,F)
\]
of objects in $\cC$ with a composition law which is associative and commutative up to coherent
homotopy. This composition law is captured in the datum of a simplicial space
$\cC_{\bullet}$ whose $k$-cells are given by the classifying space
\[
	\cC_{k} := \coprod_{n_1,\ldots,n_k} \left( \prod_{S \subseteq\{1,\ldots,k\}} \EGL(\sum_{i \in
	S} n_i, F) \right) / \GL(n_1, F) \times \GL(n_2, F) \times \cdots \times  \GL(n_k, F)
\]
of $k$-tuples of vector spaces together with choices for all possible sums. The commutativity is
reflected in additional symmetries on the simplicial space $\cC_{\bullet}$ providing it with the
structure of a $\Gamma$-space. The general theory of $\Gamma$-spaces described in
\cite{segal:categories} explicitly exhibits the ``group completion'' $\Omega |\cC_{\bullet}|$ as an
infinite loop space which, in the case at hand, yields the algebraic $K$-theory spectrum of $F$. 

Assume now that $F$ is the field $\CC$ of complex numbers. As already pointed out in
\cite{segal:categories}, we can simply modify the defining formulas for $\cC_k$, replacing the
discrete topology on the groups $\GL(n,\CC)$ with the classical Lie group topology. The group
completion of the resulting $\Gamma$-space $\overline{\cC_\bullet}$ yields the {\em topological}
complex $K$-theory spectrum $\ZZ \times \BU$.  This somewhat ad-hoc description of topological
$K$-theory as a ``topologization'' of algebraic $K$-theory can be spelled out more systematically:
namely, the classifying spaces $\cC_k$ of objects can be considered as the $\CC$-valued points of corresponding
moduli {\em stacks} $\M_k$ of objects. One can then define the complex topologization
$\M_{\CC}$ of any pre-stack $\M$ via the homotopy left Kan extension of the apparent topologization functor
\[
	\Aff \to \Top, X \mapsto X_\CC
\]
for complex affine schemes of finite type. In the above context, we have $\M_k(\CC) \simeq \cC_k$
and $(\M_k)_{\CC} \simeq \overline{\cC_k}$ so that topologization becomes an intrinsically defined
operation. Based on this idea, a theory of complex topological $K$-theory for $\CC$-linear
differential graded categories has been developed in \cite{blanc:ktheory}, based on a proposal of To\"en (see e.g. \cite{toenktop}). It yields a functor
\[
	K^{\topp}: \Cat_{\on{dg}}(\CC) \lra \on{\Sp},
\]
where $\Sp$ denotes the category of symmetric spectra.

In this work, we study the topological $K$-theory of equivariant singularity categories. In Sections~\ref{dg} -~\ref{mf}, we provide some background on dg categories, singularity categories, and graded matrix factorizations. In Sections~\ref{topktheory} and~\ref{ABS}, we establish results concerning the topological $K$-theory of dg singularity categories of weighted projective complete intersections, with special attention to weighted projective hypersurfaces.

{\bf Results.} Given a graded Gorenstein algebra $R$ over a field $\k$, let $\Ds(R)$ denote its dg singularity category (see Section~\ref{LGCY}). Let $w_1, \dots, w_n$ be positive integers, let $Q=\C[x_1, \dots, x_n]$ be the graded ring with $|x_i| = w_i$ for all $i$, let $f_1, \dots, f_c$ be a regular sequence of homogeneous elements in $Q$, and let $S=Q/(f_1, \dots, f_c)$. We prove the following (Theorem~\ref{thm:ci} below):

\begin{thm}
\label{introci}
Assume that $\sum |f_j| \le \sum w_i$. Equip the topological spaces $\cup_j\{f_j \ne 0\} \subseteq \CC^n$ and $\C^c \setminus \{0\}$ with the $S^1$ actions given by
$$z \cdot (x_1, \dots, x_n)= (z^{w_1}x_1, \dots, z^{w_n}x_n) \text{, } z \cdot (u_1, \dots, u_c) = (z^{|f_1|}u_1, \dots, z^{|f_c|}u_c).$$
 Then there is a canonical exact triangle of spectra
	\[
		K^{\topp}(\Ds(S)) \lra K_{S^1}(\CC^c\setminus \{0\})
		\overset{}{\lra} K_{S^1}(\cup_j \{f_j \ne 0\}) \overset{+1}{\lra}
	\]
\end{thm}

Here, $\K_{G}(-)$ denotes the {\em representable} $G$-equivariant complex topological $K$-theory spectrum,
as defined in \cite[Chapter XIV]{may:equivariant} (as opposed to the $K$-theory with compact support used in
\cite{segal:equivariant}).

When $c = 1$, we use the arguments in the proof of Theorem~\ref{introci}, along with applications of Bott periodicity and Kn\"orrer periodicity, to prove the following (Theorem~\ref{main} below):

\begin{thm}\label{thm:main} Let $f \in Q$ be a homogeneous polynomial of degree $d$. Fix a factorization $d = km$, where $k,m$ are positive integers. There is a canonical exact
	triangle of spectra 
	\[
			\K_{\mu_m}(\CC^n,F_f) \lra K^{\topp}(\Ds(Q/(f))) \lra K^{\topp}(\Ds(Q[u]/(f + u^k))) 
			\overset{+1}{\lra},
	\]
	where $|u| = m$, $F_f: = f^{-1}(1) \subseteq\CC^n$ denotes the affine Milnor fiber of $f$ (see Definition (1.12) in Chapter 3 of \cite{dimca2012singularities}), and the group $\mu_m \subseteq\CC^*$ of $m^{\on{th}}$ roots of unity acts on $F_f$ via monodromy.
\end{thm}

If we take $k=1$ in Theorem \ref{thm:main}, the ring $Q[u]/(f + u^k)$ is regular, and hence $\Ds(Q[u]/(f + u^k))$ is the trivial dg category. Thus, we obtain:

\begin{thm}\label{thm:ktopmf} 
	Using the notation of Theorem~\ref{thm:main}, there is a canonical weak equivalence of spectra
	\[
		 K_{\mu_d}(\CC^n, F_f) \xrightarrow{\simeq} K^{\topp}(\Ds(Q/(f))).
	\]
\end{thm}

Suppose $w_i = 1$ for all $i$, let $q_n = x_1^2 + \cdots + x_n^2 \in Q$, and let $R_n = \C[x_1, \dots, x_n]/(q_n)$. Setting $k=2$ in Theorem~\ref{thm:main}, we
obtain an exact triangle
\begin{equation}\label{eq:abs}
		K^{\topp}(\Ds(R_n)) \lra 
		K^{\topp}(\Ds(R_{n-1})) \lra
		\widetilde{K}(S^{n-1})\overset{+1}{\lra},
\end{equation}
where $\widetilde{K}$ denotes (nonconnective) reduced topological $K$-theory, and we use Kn\"orrer periodicity to reformulate the second term. Passing to $K_0$, we obtain a result
of \cite{atiyah-bott-shapiro} expressing the complex topological $K$-theory of spheres in terms of Grothendieck
groups of Clifford modules (see Section~\ref{ABS}). The triangle \eqref{eq:abs} provides a spectrum level version of this
statement which was one of the initial motivations of this work. 
We point out that the idea to relate
localization sequences like the one in Theorem \ref{thm:main} to the Atiyah-Bott-Shapiro construction has
already appeared in \cite{swan1985k}. Our contribution is the proposal to use {\em topological} (as
opposed to algebraic) $K$-theory of dg categories, which is crucial for an interpretation on the
level of spectra (as opposed to Grothendieck groups).

Our approach to studying the topological $K$-theory of matrix factorization categories is inspired by
their appearance in physics as categories of boundary conditions for so-called Landau-Ginzburg
models. The Landau-Ginzburg/Calabi-Yau (LG/CY) correspondence (see e.g. \cite{witten:phases}), or rather its
mathematical incarnation in terms of triangulated categories as established in
\cite{orlov2009derived}, allows us to relate equivariant matrix factorization categories to derived
categories of coherent sheaves. The topological $K$-theory of the latter categories can be
understood systematically via localization sequences and the comparison results of
\cite{blanc:ktheory, halpern-pomerleano:equivariant}.
\vskip\baselineskip
\noindent {\bf Acknowledgements.} We would like to thank Anthony Blanc, Ed Segal, and Bertrand To\"en for helpful discussions related to this work.

\section{Differential graded categories}
\label{dg}
Let $\k$ be a field. In this section, we recall some background on $\k$-linear differential graded categories. We will begin by introducing some basic notions, following To\"en's introduction to the subject in \cite{bertrand2011lectures}.

\subsection{Preliminaries}

A \emph{$\k$-linear differential graded (dg) category} $\fC$ consists of the following data:
\begin{itemize}
\item a set $\Ob(\fC)$ of objects;
\item for every ordered pair $(X,Y)$ of objects, a cochain complex $\Hom_\fC(X, Y)$ of $\k$-vector spaces;
\item for every ordered triple $(X, Y, Z)$ of objects, a morphism $$\Hom_\fC(X,Y) \otimes_\k \Hom_\fC(Y,Z) \to \Hom_\fC(X,Z)$$ of complexes;
\item for every object $X$, a morphism $\id_X: \k \to \Hom_\fC(X,X)$ of complexes.
\end{itemize}
This data is required to be compatible in the evident way; see Section 2.3 of \cite{bertrand2011lectures} for details.

A \emph{dg functor} $\fC_1 \to \fC_2$ consists of 
\begin{itemize}
\item a map $F: \Ob(\fC_1) \to \Ob(\fC_2)$ of sets;
\item for each ordered pair $(X,Y)$ of objects of $\fC_1$, a morphism $\Hom_{\fC_1}(X,Y) \to \Hom_{\fC_2}(F(X), F(Y))$ of complexes.
\end{itemize}
Again, we refer the reader to Section 2.3 of \cite{bertrand2011lectures} for the conditions this data must satisfy. Denote by $\Catdg(\k)$ the category of $\k$-linear dg categories with morphisms given by dg functors.

Given a dg category $\fC$, the \emph{homotopy category} of $\fC$, denoted by $[\fC]$, is defined to be the category with the same objects as $\fC$ and morphisms given by $\Hom_{[\fC]}(X,Y):=H^0\Hom_{\fC}(X,Y)$. A dg functor $F: \fC_1 \to \fC_2$ evidently induces a functor $[F]: [\fC_1] \to [\fC_2]$ on homotopy categories.

A dg functor $F$ is a \emph{quasi-equivalence} if
\begin{itemize}
\item The morphism $\Hom_{\fC_1}(X,Y) \to \Hom_{\fC_2}(F(X), F(Y))$ of complexes is a quasi-isomorphism for every ordered pair $(X,Y)$ of objects of $\fC_1$, and
\item $[F]$ is essentially surjective.
\end{itemize}



In this paper, we will often wish to consider dg categories only up to quasi-equivalence, and so we will work in the $\infty$-category $\Lqe$ of dg categories localized along quasi-equivalences. A sequence $\fC_1 \xrightarrow{} \fC_2 \xrightarrow{} \fC_3$ of dg functors is called a \emph{localization sequence} if 
\[
	\xymatrix{
		\fC_1 \ar[r]\ar[d] & \fC_2 \ar[d]\\
	0 \ar[r] & \fC_3}
\]
is a pushout square in the $\infty$-category $\Lqe$ (i.e. a homotopy pushout square with
respect to the model structure on $\Catdg(\k)$ in which the weak equivalences are quasi-equivalences; see Section 3.2 of \cite{bertrand2011lectures} for details). 

For instance, given a small $\k$-linear abelian category $\A$, denote by $\Ch^b(\A)$ the $k$-linear
dg category of bounded cochain complexes of objects in $\A$, and denote
by $\Ch^b_{\on{ac}}(\A) \subseteq\Ch^b(\A)$ the full dg subcategory spanned by the acyclic complexes. Then the \emph{dg quotient}
\[
	\D^b(\A) = \Ch^b(\A)/\Ch^b_{\on{ac}}(\A) 
\]
is defined via a pushout square
\[
	\xymatrix{
		\Ch^b_{\on{ac}}(\A) \ar[r]\ar[d] & \Ch^b(\A)\ar[d]\\
	0 \ar[r] & \D^b(\A)}
\]
in the $\infty$-category $\Lqe$. $\D^b(\A)$ is called the {\em bounded derived dg category of $\A$}. Note that $\on{D}^b(\A):=[\D^b(\A)]$ is the ordinary bounded derived category of $\A$. 

\subsection{Topological $K$-theory of dg categories}
\label{ktop}
As discussed in the introduction, Blanc has recently introduced a notion of \emph{topological $K$-theory of $\C$-linear dg categories} (\cite{blanc:ktheory}). We now recall the construction in detail. Let $\Aff$ denote the category of affine $\C$-schemes of finite type, let $\on{Set}_\Delta$ denote the category of simplicial sets, and let 
$$( - )_\C: \Aff \to \on{Set}_\Delta$$ 
denote the functor which sends a scheme $X$ to (the singular simplicial set associated to) its complex points equipped with the analytic topology. Let $h: \Aff \to \on{Fun}(\Aff^{\on{op}}, \on{Set}_\Delta)$ denote the Yoneda embedding 
$$X \mapsto (Y \mapsto \Hom(X, Y)),$$ 
where $\Hom(X, Y)$ refers to the constant simplicial set associated to the set $\Hom(X, Y)$. We denote by $\widetilde{(-)}_\C$ the $\on{Set}_\Delta$-enriched left Kan extension of $(-)_\C$ along $h$:
$$
\xymatrix{
\Aff \ar[d]^-{h} \ar[r]^-{(-)_\C} &  \on{Set}_\Delta \\
 \on{Fun}(\Aff^{\on{op}}, \on{Set}_\Delta) \ar[ru]_-{\widetilde{(-)}_\C}. \\
}
$$

As above, denote by $\Sp$ the category of symmetric spectra. Given a simplicial set $S$, let
$\Sigma^{\infty}(S)_+$ denote the suspension spectrum of $S$ with base
point adjoined, and denote the induced functor $\on{Fun}(\Aff^{\on{op}}, \on{Set}_\Delta) \to
\on{Fun}(\Aff^{\on{op}}, \Sp)$ by $\Sigma^{\infty}(-)_+$, as well. Let $T$ denote the
$\Sp$-enriched left Kan extension of $\Sigma^{\infty}(-)_+ \circ \widetilde{(-)}_\C:
\on{Fun}(\Aff^{\on{op}}, \on{Set}_\Delta) \to \Sp$ along $\Sigma^{\infty}(-)_+$:
$$
\xymatrix{
\on{Fun}(\Aff^{\on{op}}, \on{Set}_\Delta) \ar[d]^-{\Sigma^{\infty}(-)_+} \ar[rrr]^-{\Sigma^{\infty}(-)_+ \circ \widetilde{(-)}_\C} & & &  \Sp \\
\on{Fun}(\Aff^{\on{op}}, \Sp). \ar[rrru]_-{| - |} \\
}
$$
Finally, we denote by $|-|:= \on{L}T$ the left derived functor of $T$ with respect to the
\'etale $\AA^1$-local model structure on the category of spectral presheaves on $\Aff^{\on{op}}$. 

Now, let $\fC$ be a $\C$-linear dg category, and consider the object $\underline{K}(\fC)$ of $\on{Fun}(\Aff^{\on{op}}, \Sp)$ given by $\Spec(A) \mapsto K^{\alg}(A \otimes_\C \fC)$, where $K^{\alg}( - )$ denotes the nonconnective algebraic $K$-theory spectrum.

\begin{defn}[\cite{blanc:ktheory} Definition 4.1]
The \emph{semi-topological $K$-theory spectrum of $\fC$} is the spectrum $|\underline{K}(\cC)|$, denoted $K^{\on{st}}(\fC)$.
\end{defn}

\begin{rem}
This definition is inspired by Friedlander-Walker's notion of semi-topological $K$-theory for quasi-projective complex algebraic varieties (\cite{friedlander2002semi}).
\end{rem}

Let $\bf{bu}$ denote the connective topological $K$-theory spectrum of the point. By Theorem 4.6 of \cite{blanc:ktheory}, there exists a canonical isomorphism $K^{\on{st}}(\C) \cong \bf{bu}$ in the homotopy category $Ho(\Sp)$. It follows that $K^{\on{st}}(\fC)$ is a $\bf{bu}$-module. Choose a generator $\beta \in \pi_2(\bf{bu})$.

\begin{defn}[\cite{blanc:ktheory} Definition 4.13]
The \emph{topological $K$-theory spectrum of $\fC$} is defined to be the Bott-inverted semi-topological $K$-theory spectrum $K^{\on{st}}(\fC)[\beta^{-1}]$, denoted $K^{\topp}(\fC)$. We set $K_i^{\topp}(\fC):=\pi_i K^{\topp}(\fC)$. 
\end{defn}


All of the key properties of $K^{\topp}(-)$ which we will use in this paper are expressed in the following theorem:

\begin{thm}[\cite{blanc:ktheory} Theorem 1.1]
The functor $K^{\topp}: \Cat_{\on{dg}}(\CC) \lra \on{\Sp}$ enjoys the following properties:
\begin{enumerate}
	\item $K^{\topp}$ maps Morita equivalences to weak equivalences of spectra;
	\item $K^{\topp}$ maps localization sequences to exact triangles of spectra;
	\item for a separated scheme $X$ of finite type over $\Spec(\CC)$, 
		we have $K^{\topp}(\Perf(X)) \simeq K(X_\C)$. Here, $\Perf(X)$ denotes the full dg subcategory of $\D^b(\on{coh}(X))$ spanned by perfect complexes, and $K(X_\C)$ denotes the nonconnective topological $K$-theory spectrum of the space $X_\CC$ of complex points of $X$.
\end{enumerate}
\end{thm}

There is a canonical map $K^{\alg}(\fC) \to K^{\topp}(\fC)$. When $\fC = \Perf(X)$ for some separated finite type $\C$-scheme $X$, this map agrees with the usual map $K^{\alg}(X) \to K(X_\C)$. We recall that this map is typically far from being an isomorphism, as the following example illustrates:

\begin{exm}
	Let $E$ be an elliptic curve over $\C$. Then $E_\C$ is homeomorphic to a torus, so
	$K_0^{\topp}(\Perf(E)) \cong \Z^{\oplus 2}$. On the other hand, by Example 2.1.9 of
	\cite{schlichting2011higher}, $K_0^{\alg}(\Perf(E)) \cong \Z \oplus \on{Pic}(E)$, and
	$\on{Pic}(E)$ is in bijection with the set $E_\C$.
\end{exm}

\section{A theorem of Orlov}
\label{LGCY}

In this section, we discuss a result of Orlov (\cite[Theorem 2.5]{orlov2009derived}) which we will use in the proof of Theorem \ref{thm:ci}. Let $\k$ be a field, and let $R = \oplus_{i \ge 0} R_i$ be a Noetherian graded $\k$-algebra such that $R_0 = \k$. We assume that $R$ is Gorenstein; that is, 
\begin{itemize}
	\item $R$ has finite injective dimension $n$ over itself,
	\item $\RHom_R(\k,R) \cong \k(a)[-n]$ for some $a \in \ZZ$ which we call the {\em Gorenstein parameter}. Here, the round brackets denote a shift in the internal grading, and the square brackets denote a shift in the cohomological grading.
	\end{itemize}
We denote by $\m_R = R_{\ge 1}$ the irrelevant ideal.

\subsection{Noncommutative projective geometry}

Consider the abelian category $\gr R$ of finitely generated graded right $R$-modules 
and the full subcategory $\tor R \subseteq\gr R$ consisting of modules which are finite
dimensional over $\k$. We introduce the quotient category 
\[
	\qgr R = \gr R / \tor R
\]
and obtain, as shown in \cite{artin-zhang:noncommutative}, an adjunction
\[
	p: \gr R \llra \qgr R: \Gamma,
\]
where $p$ is the projection functor and the right adjoint $\Gamma$ is fully faithful. We further
consider the adjunction
\[
	\iota: \grn R \llra \gr R: \tau_{\ge 0},
\]
where $\iota$ is the inclusion of the full subcategory $\grn R$ consisting of modules $M$ satisfying
$M_i = 0$ for $i <0$, and $\tau_{\ge 0}$ is the degree truncation functor. Composing the two
adjunctions yields the adjunction
\[
	p \iota: \grn R \llra \qgr R: \Gamma_{\ge 0},
\]
where $\Gamma_{\ge 0}$ is still fully faithful.

\begin{prop} 
\label{square}
Let $\varphi: R \to S$ be a homomorphism of graded Gorenstein algebras such that $S$ is a finite
	$R$-module of finite projective dimension. Then there is an induced commutative square 
	\begin{equation}\label{eq:square-qgr}
		\xymatrix{
			\D^b(\qgr S) \ar[d]_{\R\Gamma_S} \ar[r]^{\overline{\varphi^*}}
			&\ar[d]^{\R\Gamma_R} \D^b(\qgr R)\\
			\D^b(\gr S) \ar[r]^{\varphi^*}  &  \D^b(\gr R) 
		}
	\end{equation}
	in the $\infty$-category $\Lqe$.
\end{prop}
\begin{proof}
	We first produce the square on the level of homotopy categories. Namely, due to the
	assumptions on $\varphi$, the functor
	\[
		\L\varphi_! = - \otimes^L_R S: \Db(\gr R) \lra \Db(\gr S)
	\]
	preserves complexes which are quasi-isomorphic to complexes of torsion modules and hence induces a
	square
	\[
		\xymatrix{
			\Db(\qgr S) & \ar[l]_{\overline{\L\varphi_!}} \Db(\qgr R)\\
			\Db(\gr S) \ar[u]^{p} & \ar[l]_{\L \varphi_!} \Db(\gr R)\ar[u]_{p} 
		}
	\]
	which commutes up to natural isomorphism. Passing to right adjoints, we obtain the commutative
	square 
	\begin{equation}\label{eq:homotopy-square}
			\xymatrix{
				\Db(\qgr S) \ar[d]_{\R\Gamma_S} \ar[r]^{\overline{\varphi^*} } &\ar[d]^{\R\Gamma_R}
				\Db(\qgr R)\\
				\Db(\gr S) \ar[r]^{\varphi^*}  &  \Db(\gr R). \\
			}
	\end{equation}
	To lift this to a square in $\Lqe$, note that restricting the functor 
	\[
		p: \D^b(\gr S) \lra \D^b(\qgr S)
	\]
	to the full dg subcategory spanned by the objects in $\R \Gamma_S(\Db(\qgr S))$ yields a
	quasi-equivalence. Therefore, in $\Lqe$, we may identify the dg category $\D^b(\qgr S)$ with the full
	dg subcategory in $\D^b(\gr S)$ spanned by $\R \Gamma_S(\Db(\qgr S))$. By
	\eqref{eq:homotopy-square}, the dg functor $\varphi^*$ maps this subcategory to the corresponding
	subcategory in $\D^b(\gr R)$. This yields the desired square \eqref{eq:square-qgr} in $\Lqe$.
\end{proof}

In the setting of Proposition~\ref{square}, the truncation functors $\tau_{\ge 0}$ induce a square 
\[
	\xymatrix{
		\D^b(\gr S) \ar[d]_{\tau_{\ge 0}} \ar[r]^{\varphi^*} &\ar[d]^{\tau_{\ge 0}} \D^b(\gr
		R)\\
		\D^b(\grn S) \ar[r]^{\varphi^*}  &  \D^b(\grn R) 
	}
\]
which we may compose with \eqref{eq:square-qgr} to obtain the commutative square 
\begin{equation}\label{eq:square-qgrn}
		\xymatrix{
			\D^b(\qgr S) \ar[d]_{\R\Gamma_{\ge 0}} \ar[r]^{\overline{\varphi^*}}
			&\ar[d]^{\R\Gamma_{\ge 0}} \D^b(\qgr R)\\
			\D^b(\grn S) \ar[r]^{\varphi^*}  &  \D^b(\grn R) 
		}
\end{equation}
in $\Lqe$.

\subsection{Graded singularities}
\label{orlovsection}
We denote by $\Perf R \subseteq\D^b(\gr R)$ the full dg subcategory spanned by those objects
which are isomorphic in $\on{D}^b(\gr R)$ to bounded complexes of finite projective $R$-modules. We
introduce the dg quotient
\[
	\Ds(R) := \D^b(\gr R) / \Perf R,
\]
the \emph{graded singularity category} of $R$. Set $\on{D}^{\sg}(R):=[\Ds(R)]$.
\begin{prop} Let $\varphi: R \to S$ be a homomorphism of graded Gorenstein algebras such that $S$ is a finite
	$R$-module of finite projective dimension. Then there is an induced commutative square 
	\begin{equation}\label{eq:square-sing}
		\xymatrix{
			\D^b(\gr S) \ar[d]_-{\pi} \ar[r]^-{\varphi^*} & \D^b(\gr R)\ar[d]^-{\pi}\\
			\Ds(S) \ar[r]^-{\overline{\varphi^*}}  &  \Ds(R). 
		}
	\end{equation}
	in the $\infty$-category $\Lqe$.
\end{prop}
\begin{proof}
	Clear since, due to the assumptions on $\varphi$, the dg functor $\varphi^*: \D^b(\gr S) \to
	\D^b(\gr R)$ preserves perfect complexes.
\end{proof}

We denote by $\Lambda_R$ the composite of the dg functors
\[
	\D^b(\qgr R) \overset{\R \Gamma_{\ge 0}}{\lra} \D^b(\grn R) \overset{\pi \circ
	\iota}{\lra} \Ds(R).
\]
The following is a celebrated theorem of Orlov:

\begin{thm}[\cite{orlov2009derived} Theorem 2.5]
\label{orlov} Let $R$ be a graded Gorenstein algebra with Gorenstein parameter $a$. Then 
	\begin{enumerate}
		\item if $a \ge 0$, $\Lambda_R$ is a localization with
			kernel given by the full dg subcategory of $\D^b(\qgr R)$ spanned by the
			exceptional collection $\langle p R(-a+1), \cdots, p R \rangle$,
		\item if $a \le 0$, $\Lambda_R$ is fully faithful with orthogonal complement in
			$\Ds(R)$ given by the exceptional collection $\langle \pi k(-a+1), \cdots, \pi
			k \rangle$.
	\end{enumerate}
\end{thm}

Let $\varphi: R \to S$ be a homomorphism of graded Gorenstein algebras such that $S$ is a finite
$R$-module of finite projective dimension. 
Composing \eqref{eq:square-sing} and \eqref{eq:square-qgrn}, we arrive at the fundamental commutative square
\begin{equation}\label{eq:fund-square}
		\xymatrix{
			\D^b(\qgr S) \ar[d]_{\Lambda_S} \ar[r]^{\overline{\varphi^*}} & \D^b(\qgr R)\ar[d]^{\Lambda_R}\\
			\Ds(S) \ar[r]^{\overline{\varphi^*}}  &  \Ds(R) }
\end{equation}
in $\Lqe$.

\section{$\Z$-graded matrix factorizations}
\label{mf}
Let $\k$ be a field, and let $Q$ denote the graded polynomial ring $\k[x_1, \dots, x_n]$, where $|x_i|$ is a positive integer for each $i$. Let $f \in Q$ be nonzero and homogeneous of degree $d$, and let $S$ denote the graded hypersurface ring $Q/(f)$.

It is often convenient to use an alternative model for the graded singularity category $\Ds(S)$: the dg category of graded matrix factorizations of $f$. In this section, we introduce graded matrix factorizations, and we briefly discuss some of their key properties.

\begin{defn}
\label{mfs}
The $\k$-linear dg category $\MF^\Z(f)$ of \emph{graded matrix factorizations} of $f$ over $Q$ is given by the following:
\begin{itemize}
\item An object of $\MF^\Z(f)$ is a pair of graded free $Q$-modules $F_0, F_1$ equipped with a degree $d$ map
$$s_0: F_0 \to F_1$$
and a degree $0$ map
$$s_1: F_1 \to  F_0$$
satisfying $s_0 s_1 = f \cdot \id_{F_1}$ and $s_1 s_0 = f \cdot \id_{F_0}$. We will denote objects by $F_0 \darrow{s_0}{s_1} F_1$.

\item The morphism complex $\Hom_{\MF}((F_0 \darrow{s_0}{s_1} F_1) , (F_0' \darrow{s_0'}{s_1'} F_1'))$ has underlying graded $\k$ - vector space given by
$$\Hom_{\gr Q}(F_0, F_0'(ld)) \oplus \Hom_{\gr Q}(F_1, F_1'(ld))$$
in degree $2l$, and
$$\Hom_{\gr Q}(F_0, F_1'((l+1)d)) \oplus \Hom_{\gr Q}(F_1, F_0'(ld))$$
in degree $2l + 1$. 

The differential $\partial$ on $\Hom_{\MF}((F_0 \darrow{s_0}{s_1} F_1) , (F_0' \darrow{s_0'}{s_1'} F_1'))$ is given by 
$$\partial(\alpha) = (s_0' + s_1') \alpha - (-1)^{|\alpha|} \alpha (s_0 + s_1).$$ 

\end{itemize}
\end{defn}

\begin{rem}
The data of a graded matrix factorization of $f$ is equivalent to a sequence of degree $0$ maps
$$\cdots \to F_0(-2d) \xrightarrow{s_0} F_1(-d) \xrightarrow{s_1} F_0(-d) \xrightarrow{s_0} F_1 \xrightarrow{s_1} F_0 \xrightarrow{s_0} F_1(d) \to \cdots$$
of graded free $Q$-modules such that consecutive maps are given by multiplication by $f$. Such sequences are obviously not complexes, but thinking of them as such can yield useful intuition; for instance, observe that a degree $m$ morphism of graded matrix factorizations that is a cocycle is, from this point of view, precisely a degree $m$ map of ``complexes".
\end{rem}

\begin{rem}
Since $f$ is a non-zero-divisor in $Q$, the ranks of the free $Q$-modules $F_1$ and $F_0$ in Definition~\ref{mfs} are equal.
\end{rem}

There is an equivalence of categories
$$\Phi: [\MF^{\Z}(f)] \xrightarrow{\cong} \on{D}^{\sg}(S)$$
given, on objects, by $(F_0 \darrow{s_0}{s_1} F_1) \mapsto \coker(s_1)$, where the latter is considered as a complex concentrated in degree 0. This theorem is essentially due to work of  Buchweitz and Eisenbud in \cite{buchweitz1986maximal} and \cite{eisenbud1980homological}; it is proven explicitly by Orlov in Theorem 3.10 of \cite{orlov2009derived}. Note that $\Phi$ lifts to a quasi-equivalence $\MF^\Z(f) \xrightarrow{\simeq} \Ds(S)$ of dg categories.

In particular, $[\MF^\Z(f)]$ has a canonical triangulated structure. The shift functor applied to a matrix factorization $F=(F_0 \darrow{s_0}{s_1} F_1)$ is given by 
$$F[1]:= (F_1(d) \darrow{-s_1}{-s_0} F_0).$$
One may also apply a grading shift to a matrix factorization $F=(F_0 \darrow{s_0}{s_1} F_1)$; if $l \in \Z$, define
$$F(l):= (F_0(l) \darrow{s_0}{s_1} F_1(l)).$$
Observe that $F[2] = F(d)$.

\begin{exm}[\cite{ballard2014category} page 28]
\label{power}
Suppose $n=1$ and $f = x_1^d$ (so $w_1 = 1$). Then there exists a quasi-equivalence
$$\MF^\Z(f) \xrightarrow{\simeq} \D^b(A_{d-1}),$$
where $\D^b(A_{d-1})$ denotes the dg bounded derived category of the $A_{d-1}$-quiver. In particular, if $d = 2$, one has a quasi-equivalence 
$$\MF^\Z(f) \xrightarrow{\simeq} \Perf (\k).$$
The matrix factorization $\k[x_1] \darrow{x_1}{x_1} \k[x_1](-1)$ may be taken to correspond to the complex in $\Perf(\k)$ with $\k$ concentrated in degree 0. 
\end{exm}

\begin{exm}
\label{uv}
Suppose $n =2$ and $f = x_1x_2$, where $w_1 = 1 = w_2$. By Theorem~\ref{orlov}, there exist isomorphisms $\MF^\Z(f) \cong \D^b(\Proj(S)) \cong \Perf (\k) \times \Perf (\k)$ in $\Lqe$. The objects
$$\k[x_1, x_2] \darrow{x_1}{x_2} \k[x_1, x_2](-1) \text{, } \k[x_1, x_2] \darrow{x_2}{x_1} \k[x_1, x_2](-1)$$
of $[\MF^\Z(f)]$ may be taken to correspond to the two copies of the complex with $\k$ concentrated in degree 0 in $\Perf(\k) \times \Perf(\k)$.
\end{exm}

Let $Q':=\k[y_1, \dots, y_{m}]$, where each $y_i$ has some positive integer weight, and suppose $f'$ is a degree $d$ element of $Q'$. If $F:=(F_0 \darrow{s_0}{s_1} F_1)$ and $F':=(F_0' \darrow{s_0'}{s_1'} F_1')$ are objects of $\MF^\Z(f)$ and $\MF^\Z(f')$, we define their \emph{tensor product}, which is an object in $\MF^\Z(f + f')$, as follows:

$$F \otimes_{\MF} F' :=(F_0 \otimes F_0') \oplus (F_1 \otimes F_1')(d) 
\darrow{
\begin{pmatrix}  
\id \otimes s_0' & s_1 \otimes \id \\ 
-s_0 \otimes \id & \id \otimes s_1'
\end{pmatrix}
}
{
\begin{pmatrix}  
\id \otimes s_1' & -s_1 \otimes \id \\ 
s_0 \otimes \id & \id \otimes s_0'
\end{pmatrix}
} 
(F_0 \otimes F_1') \oplus (F_1 \otimes F_0').$$

The tensor product of matrix factorizations may be viewed as a dg functor. To make this precise, we recall the notion of a tensor product of dg categories. Given two $\k$-linear dg categories $\fC_1$, $\fC_2$, we define the \emph{tensor product} $\fC_1 \otimes_\k \fC_2$ to be the dg category with objects $\Ob(\fC_1) \times \Ob(\fC_2)$, morphism complexes 
$$\Hom_{\fC_1 \otimes_\k \fC_2}((X_1, X_2), (Y_1, Y_2)):= \Hom_{\fC_1}(X_1, Y_1) \otimes_\k \Hom_{\fC_2}(X_2, Y_2),$$
and the evident composition and units. It is easy to check that $-\otimes_{\MF} -$ yields a dg functor
$$\MF^\Z(f) \otimes_\k \MF^\Z(f') \to \MF^\Z(f + f').$$


\begin{rem}
We point out that, unlike the tensor product functor in the setting of non-graded matrix factorizations of isolated hypersurface singularities, the map $\MF^\Z(f) \otimes_\k \MF^\Z(f') \to \MF^\Z(f + f')$ described above is not necessarily a Morita equivalence (see Section 4.4 of \cite{bertrand2011lectures} for the definition of a Morita equivalence of dg categories). For instance, take $\k = \C$, $n = 1 = m$, $d = 2$, $f = x_1^2$, and $f' = y_1^2$; Examples~\ref{power} and~\ref{uv} imply that the tensor product functor is not a Morita equivalence in this case (cf. \cite{brown2016knorrer} Proposition 2.14).
\end{rem}

We recall Kn\"orrer's periodicity theorem for non-graded matrix factorizations. Given a commutative ring $T$ and an element $h$ of $T$, let $\MF(T,h)$ denote the differential $\Z/2\Z$-graded category of matrix factorizations of $h$ over $T$; for the definition of this category, see \cite{dyckerhoff2011compact} Definition 2.1. Let $[\MF(T,h)]$ denote the homotopy category of $\MF(T,h)$.

\begin{thm}[\cite{knorrer1987cohen} Theorem 3.1]
\label{kpint}
If $k$ is algebraically closed, $\charr(k) \ne 2$, and $h \in (x_1, \dots, x_n)
\subseteq \k[[x_1, \dots, x_n]]$, there is an equivalence of categories
$$[\MF(\k[[x_1, \dots, x_n]],h)] \xrightarrow{\cong}
[\MF(\k[[x_1, \dots, x_n,u,v]],h + u^2 + v^2)].$$
\end{thm}

The following is a version of Kn\"orrer periodicity for graded matrix factorizations:

\begin{thm}
\label{gradedknorrer}
Let $l$ be an integer such that $1 \le l \le d$, and let $\k[u,v]$ be the graded polynomial ring with $|u| = l$ and $|v| = d-l$. Let $X$ be the object $\k[u,v] \darrow{u}{v} \k[u,v](l-d) $ of $\MF^\Z(uv)$. Then the dg functor
$$K: \MF^\Z(f) \to \MF^\Z(f + uv)$$
given, on objects, by $- \otimes_{\MF} X$ and, on morphisms, by $- \otimes \id_X$
is a quasi-equivalence.
\end{thm}

This result appears to be well-known (see, for instance, Remark 2.9 of \cite{orlov2006triangulated}); nevertheless, we provide a proof in Appendix~\ref{kappendix}.


\begin{rem}
Suppose $\charr(\k) \ne 2$, $\k$ contains a square root of $-1$, $d$ is even, and $|u| = \frac{d}{2} = |v|$. Then the dg functor
$$- \otimes_{\MF} (\k[u,v] \darrow{u + iv}{u - iv} \k[u,v](-\frac{d}{2})): \MF^\Z(f) \to \MF^\Z(f + u^2 + v^2)$$
is evidently also a quasi-equivalence. We will denote this functor by $K$ as well.
\end{rem}

\begin{rem}
One may replace $X$ with its shift or any of its grading twists in the statement of Theorem~\ref{gradedknorrer}, and the result still holds.
\end{rem}
\section{Topological $K$-theory of equivariant singularity categories}
\label{topktheory}

Let $w_1, \dots, w_n$ be positive integers, and let $Q = \CC[x_1,\dots,x_n]$ be the graded polynomial ring such that $|x_i| = w_i$. Let $f \in Q$ be
a quasihomogeneous polynomial of degree $d$. The map 
\[
	f: \CC^n \setminus \{f^{-1}(0)\} \lra \CC^*
\]
is a fibration which we call the {\em affine Milnor fibration}, following Definition (1.12) in Chapter 3 of \cite{dimca2012singularities}. Due to the quasihomogeneity of $f$,
for every $\lambda \in \CC^*$ and $a \in \CC^*$, scalar multiplication by $\lambda$ provides a
diffeomorphism $f^{-1}(a) \xrightarrow{\cong} f^{-1}(\lambda^d a)$. It is straightforward to see that
multiplication by $\exp(2\pi i/d)$ induces the monodromy operator on each fiber $f^{-1}(a)$. 
This operator generates an action of the group $\mu_d \subseteq\CC^*$ of $d^{\on{th}}$ roots of unity which
we call the {\em monodromy action}. We refer to the fiber $F_f := f^{-1}(1)$ as the {\em affine
Milnor fiber}. The affine Milnor fibration is homotopy equivalent to the Milnor fibration of $f$ as 
defined in \cite{milnor1968singular} (see \cite{dimca2012singularities}). 

One of our goals is to construct a canonical weak equivalence from the relative $\mu_d$-equivariant topological $K$-theory spectrum of the pair $(\C^n, F_f)$ to the topological $K$-theory spectrum of the dg category $\Ds(Q/(f))$. We will first establish a more
general result for complete intersections which holds whenever the weighted number of variables is
large enough. We will then use Kn\"orrer periodicity to eliminate these conditions in the case of
hypersurface singularities.

\subsection{Complete intersections}
\label{sec:complete}

Consider a
regular sequence of quasihomogeneous polynomials $f_1, \ldots, f_c$ of degrees $|f_j| = d_j$. 
For each $j$, choose a factorization $d_j = k_j m_j$, where $k_j, m_j$ are positive integers, and define $g_j := f_j + u_j^{k_j} \in Q[u_1, \dots, u_c]$, where $|u_j| = m_j$. Notice that the degree of $g_j$ is again $d_j$, and $g_1, \dots, g_c$ is a regular sequence in $Q[u_1, \dots, u_c]$.

One has a morphism of Artin stacks
\[
	\underline{u}: [\CC^{n+c}/\CC^*] \lra [\CC^c/\CC^*] 
\]
given by $(x_1, \dots, x_{n}, u_1, \dots, u_c) \mapsto (u_1, \dots, u_c)$, where the $\CC^*$ actions on $\CC^{n+c}$ and $\CC^c$ are given by the weights of the variables $(x_1, \dots, x_n, u_1, \dots, u_c)$ and $(u_1, \dots, u_c)$, respectively.

We fix the following notation:
\begin{itemize}
\item $\PP^c := [\C^c \setminus \{0 \}/ \C^*]$ and $\PP^{n+c} := [\C^{n+c} \setminus \{0 \}/ \C^*]$,
\item $\sX:= [ \{g_1, \dots, g_c = 0\}/\C^*] \subseteq \PP^{n+c}$,
\item $\sZ:=  [ \{f_1, \dots, f_c = 0\}/\C^*] \subseteq \sX$,
\item $\sU := \sX \setminus \sZ $.
\end{itemize}

We have a diagram of stacks
\begin{equation}\label{eq:relmil}
		\xymatrix{ \sZ \ar[dd]^p\ar[r]^-i
		& \sX \ar@/^-7ex/[dd]_q\ar[d] & \ar[l] \sU \ar[dd]^{\Psi}\\
		& [\CC^{n+c}/\CC^*]\ar[d]^{\underline{u}} &\\ 
		[\{0\}/\CC^*] \ar[r] & [\CC^c/\CC^*] & \ar[l] \PP(\CC^c).}
\end{equation}
Note that the two squares are pullback squares, and the horizontal arrows correspond to open-closed
decompositions. The morphism $\Psi$ should be thought of as a
relative equivariant affine Milnor fiber: for every $1 \le j \le c$, there is a pullback square
\[
	\xymatrix{
		[F_{f_j}/\mu_{d_j}] \ar[r] \ar[d] & \sU \ar[d]^{\Psi}\\
		[\ast/\mu_{d_j}] \ar[r] & \PP(\C^c)
	}
\]
where $F_{f_j} := f_j^{-1}( 1)$ is the affine Milnor fiber of $f_j$ equipped with the $\mu_{d_j}$
action induced by monodromy.

Set $S=Q/(f_1, \dots, f_c)$ and $R= Q[u_1, \dots, u_c]/(g_1, \dots, g_c)$, and let $\varphi: R \to S$ be the map given by setting the $u_j$ to 0. Let $\D^b(\qgr S) \xrightarrow{\overline{\varphi^*}} \D^b(\qgr R)$ 
be the dg functor defined in Section~\ref{LGCY}. One evidently has a commutative square 

\[
		\xymatrix{
		\D^b(\qgr S) \ar[d]_{\widetilde{( - )}} \ar[r]^-{\overline{\varphi^*}}  &  \D^b(\qgr R) \ar[d]_{\widetilde{( - )}} \\
			\D^b(\coh \sZ) \ar[r]^-{i_*} & \D^b(\coh \sX) }			
\]
in $\Lqec$; it follows from Proposition 2.17 of \cite{orlov2009derived} that $\qgr R \xrightarrow{\widetilde{(-)}} \coh \sX$ is an equivalence of categories, and similarly for $S$ and $\sZ$.

Assume now that $\sum d_j \le \sum w_i$. Applying Theorem~\ref{orlov} to the square~(\ref{eq:fund-square}), we have the following commutative diagram
\[
		\xymatrix{
			\langle \OO(-a), \dots, \OO \rangle \ar[d] \ar[r]^-{} & \langle \OO(-b) , \dots, \OO \rangle \ar[d] \\
			\D^b(\coh \sZ) \ar[d]_-{\Lambda_S \circ \widetilde{( - )}^{-1}} \ar[r]^-{i_*} & \D^b(\coh \sX)\ar[d]^-{\Lambda_R \circ \widetilde{( - )}^{-1}}\\
			\Ds(S) \ar[r]^-{\overline{\varphi^*}}  &  \Ds(R) }
\]
in $\Lqec$, where $a := (\sum_{i=1}^n w_i) - (\sum_{j=1}^c d_j) - 1$, $b:= a + \sum_{j=1}^c m_j$, the top vertical maps are induced by inclusions, and the two columns are localization sequences.

Using \eqref{eq:relmil}, we may refine the top square to obtain
	\begin{equation}\label{refinement}
		\xymatrix{
			\langle \OO(-a),\cdots,\OO \rangle \ar[d]\ar[r] & \langle \OO(-b),\cdots,\OO
			\rangle \ar[d] \\ 
			\D^b(\coh [\{0\}/\CC^*]) \ar[r]\ar[d]^{p^*} & \D^b(\coh [\CC^c/\CC^*])
			\ar[d]^{q^*} \\
			\D^b(\coh \sZ))  \ar[r]^{i_*} & \D^b(\coh \sX ).
		}
	\end{equation}
	This statement is obtained by verifying it for the respective abelian categories of coherent sheaves
	and then passing to derived dg categories. Passing to topological $K$-theory and forming cofibers of all horizontal maps, we obtain the
	diagram of spectra
	\begin{equation}\label{eq:relmildiag}
			\xymatrix{
				K^{\topp}(\langle \OO(-a),\dots,\OO \rangle) \ar[d]\ar[r] & K^{\topp}(\langle \OO(-b),\cdots,\OO
				\rangle) \ar[d] \ar[r] & V \ar[d] \\ 
				K^{\topp}(\D^b(\coh [\{0\}/\CC^*])) \ar[r]\ar[d]^{p^*} & K^{\topp}(\D^b(\coh
				[\CC^c/\CC^*]))
				\ar[d]^{q^*} \ar[r] & K_{S^1}(\CC^c\setminus \{0\}) \ar[d]^-{\Psi^*}\\
				K^{\topp}(\D^b(\coh \sZ))  \ar[r]^{i_*} \ar[d]^-{\Lambda_S \circ \widetilde{( - )}^{-1}} & K^{\topp}(\D^b(\coh \sX)) \ar[r] \ar[d]^-{\Lambda_R \circ \widetilde{( - )}^{-1}} & K_{S^1}(\cup_i \{f_i \ne 0\}) \ar[d] \\
				K^{\topp}(\Ds(S)) \ar[r]^-{\overline{\varphi^*}}  &  K^{\topp}(\Ds(R)) \ar[r] & W 
			}
	\end{equation}
	where we use (the proof of) Lemma 3.6 and Theorem 3.9 of \cite{halpern-pomerleano:equivariant} to determine the middle two horizontal cofibers. Here, $K_{G}(-)$ denotes the $G$-equivariant complex topological $K$-theory spectrum as defined in \cite[Chapter XIV]{may:equivariant}, as discussed in the introduction.

\begin{prop}
\label{we}
The map $V \to K_{S^1}(\C^c \setminus 0)$ in (\ref{eq:relmildiag}) is a weak equivalence. In particular, one has an exact triangle of spectra
$$K_{S^1}(\C^c \setminus \{0\}) \xrightarrow{\Psi^*} K_{S^1}(\cup_i\{f_i \ne 0 \}) \to W \xrightarrow{+1}.$$
\end{prop}

\begin{proof}
By the properties of the functor $K^{\topp}$ discussed in Section~\ref{ktop}, if $\fC$ is generated by an exceptional collection, we have:
\begin{itemize}
\item $K_1^{\topp}(\fC) = 0$,
\item the canonical map $K_0^{\alg}(\fC) \to K_0^{\topp}(\fC)$ is an isomorphism.
\end{itemize}

Thus, since each of the categories in the top two rows of~(\ref{refinement}) is generated by an exceptional collection, it suffices to show that the induced maps
$$\ker(K_0^{\alg}(\alpha)) \to \ker(K_0^{\alg}(\beta))$$
$$\coker(K_0^{\alg}(\alpha)) \to \coker(K_0^{\alg}(\beta))$$
are isomorphisms, where $\alpha$ and $\beta$ denote the top and middle horizontal functors in~(\ref{refinement}), respectively.

Let $\on{Kos}(\underline{u}) \in \langle \OO(-b), \dots, \OO \rangle \subseteq \D^b(\coh \sX)$ denote the Koszul complex on the sequence $u_1, \dots, u_c$. Notice that $\alpha(\OO(i)) = \on{Kos}(\underline{u})(i)$ for all $-a \le i \le 0$. Similarly, $\beta(\OO(i))$ is the $i$th twist of the Koszul complex on the coordinates of $\C^c$ for all $i \in \Z$. It follows easily that $K_0^{\alg}(\alpha)$ and $K_0^{\alg}(\beta)$ are injective, and a routine calculation shows the induced map on cokernels is an isomorphism as well.
\end{proof}

We now easily deduce the following:

\begin{thm}\label{thm:ci} Assume that $\sum d_j \le \sum w_i$. Equip the topological spaces $\cup_j\{f_j \ne 0\} \subseteq \C^n$ and $\C^c \setminus \{0\}$ with the $S^1$ actions given by
$$z \cdot (x_1, \dots, x_n)= (z^{w_1}x_1, \dots, z^{w_n}x_n) \text{, } z \cdot (u_1, \dots, u_c) = (z^{d_1}u_1, \dots, z^{d_c}u_c).$$
 Then there is a canonical exact triangle of spectra
	\[
		K^{\topp}(\Ds(S)) \lra K_{S^1}(\CC^c\setminus \{0\})
		\overset{\Psi^*}{\lra} K_{S^1}(\cup_j \{f_j \ne 0\}) \overset{+1}{\lra}
	\]
\end{thm}

\begin{proof}
Take $k_j = 1$ for all $j$ (so $m_j = d_j$ for all $j$), and consider~(\ref{eq:relmildiag}). Notice that $R$ is regular in this case, and so $\Ds(R)$ is the trivial dg category. Now use Proposition~\ref{we}.
\end{proof}


\subsection{Hypersurfaces}
We now consider the case where $c=1$. Let $f \in Q$ be homogeneous of degree $d$, and set $S = Q/(f)$. Suppose $d = km$, where $k, m$ are positive integers, and set $R = Q[u]/(f + u^k),$ where $|u| = m$.

Recall that $\mu_d$ acts on the affine Milnor fiber $F_f$ via monodromy. Let $K_{\mu_d}(\C^n, F_f)$ denote the relative $K$-theory spectrum of the inclusion $F_f \into \C^n$ of $\mu_d$-spaces (i.e. the fiber of the pullback map $K_{\mu_d}(\C^n) \to K_{\mu_d}(F_f)$ of spectra). The following is immediate from~(\ref{eq:relmildiag}) and Proposition~\ref{we}:

\begin{prop} 
\label{hs}
If $\sum w_i \ge d$, there exists a canonical exact triangle of spectra
$$K_{\mu_m}(\C^n, F_f) \to K^{\topp}(\Ds(S)) \xrightarrow{\overline{\varphi^*}} K^{\topp}(\Ds(R))  \xrightarrow{+1}.$$
In particular, there exists a weak equivalence
$$K_{\mu_d}(\C^n, F_f) \xrightarrow{\simeq} K^{\topp}(\Ds(S)).$$
\end{prop}

We now wish to use Kn\"orrer periodicity (Theorem~\ref{gradedknorrer}) to eliminate the assumption $\sum w_i \ge d$ from Proposition~\ref{hs}. As above, let $\varphi: R \to S$ be given by setting $u = 0$, and define $T$ to be the dg functor
$$- \otimes_{\MF} (\C[u] \darrow{u^{k-1}}{u} \C[u](-m)): \MF^\Z(f) \to \MF^\Z(f +  u^k).$$
We will need the following:
\begin{lem}
\label{pushtensor}
One has a commutative diagram
\[
	\xymatrixcolsep{1.1in}
	 \xymatrix{
	 \MF^\Z(f) \ar[d]^-{\Phi} \ar[r]^-{T} & \MF^\Z(f + u^k) \ar[d]^-{\Phi} \\
	\Ds(S)  \ar[r]^-{\overline{\varphi^*}}  &  \Ds(R)  }	
\]
in $\Lqec$, where $\overline{\varphi^*}$ is as in Section~\ref{LGCY}, and the vertical maps $\Phi$ are the quasi-equivalences discussed in Section~\ref{mf}.
\end{lem}

\begin{proof}
Let $F:=(F_0 \darrow{s_0}{s_1} F_1)$ be an object of $\MF^\Z(f)$. Write $M = \coker(s_1)$, and let $\on{res}(M)$ denote the $R$-module obtained by restriction of scalars, so that $(\overline{\varphi^*} \circ \Phi)(F)$ is the complex with $\on{res}(M)$ concentrated in degree 0. Let $\overline{F_i} := F_i \otimes_Q R$ for $i = 0, 1$. $\on{res}(M)$ has an $R$-free resolution
$$
\cdots \rightarrow \overline{F_0}(-d-m) \oplus \overline{F_1}(-d)  \xrightarrow{\begin{pmatrix} 
u & s_1  \\
-s_0 & u^{k-1} \\
\end{pmatrix}} \overline{F_0}(-d) \oplus \overline{F_1}(-m) \xrightarrow{\begin{pmatrix} 
u^{k-1} & -s_1  \\
s_0 & u \\
\end{pmatrix}} \overline{F_0}(-m) \oplus \overline{F_1} \xrightarrow{\begin{pmatrix} u & s_1 \end{pmatrix}} \overline{F_0}.
$$

On the other hand, setting $F_i'= F_i \otimes_\C \C[u]$ for $i = 0,1$, we have
$$\Phi (F \otimes_{\MF} (\C[u](-m) \darrow{u}{u^{k-1}} \C[u](-d))) = \coker (F_0'(-d) \oplus F_1'(-m)  \xrightarrow{\begin{pmatrix} 
u^{k-1} & -s_1  \\
s_0 & u \\
\end{pmatrix}}  F_0'(-m) \oplus F_1'),$$
so that, by (a $\Z$-graded version of) Lemma 2.2.2 of \cite{buchweitz1986maximal}, there is an isomorphism 
$$\on{res}(M)[-1] \cong \Phi (F \otimes_{\MF} (Q[u](-m) \darrow{u}{u^{k-1}} Q[u](-d))) = \Phi(T(F)[-1])$$ 
in $\on{D}^{\sg}(R)$. \end{proof}

We will also need the following classical results:

\begin{thm}[Equivariant Bott Periodicity] Let $G$ be a compact Lie group and $X$ a compact $G$-space. If $V$ is a complex representation of $G$, there exists a weak equivalence
$$K_G(X) \xrightarrow{\simeq} K_G(\Sigma^V X).$$
\end{thm}

\begin{thm}[The Sebastiani-Thom Theorem] 
\label{ST}
Let $f \in \C[y_1, \dots, y_l]$ and $f' \in \C[y_1', \dots, y_l']$ be quasihomogeneous polynomials. Then there exists a homotopy equivalence
$$F_f * F_{f'} \xrightarrow{\simeq} F_{f + f'},$$
where the left side is the topological join of $F_f$ and $F_{f'}$. Further, this homotopy equivalence preserves monodromy operators.
\end{thm}

Theorem~\ref{ST} was originally proved for isolated hypersurface singularities by Sebastiani-Thom in \cite{sebastiani1971resultat} and was generalized to arbitrary quasihomogeneous polynomials by Oka in \cite{oka1973homotopy}.

Now, let $a$ and $b$ be positive integers such that $a + b = d$. Consider the graded ring $\C[v, w]$, where $|v| = a$ and $|w| = b$, so that $vw$ is homogeneous of degree $d$. Recall that the affine Milnor fiber $F_{vw}$ is equipped with the $\mu_d$-action given by $\zeta \cdot (v, w) = (\zeta^a v, \zeta^b w)$, where $\zeta = \exp(2 \pi i / d)$. One easily checks that $F_{vw}$ equivariantly deformation retracts onto the subspace $\{(x, \frac{1}{x}) \in F_{vw} \text{ : } |x| = 1 \}$. Moreover, this subspace is equivariantly homeomorphic to $S^1$ equipped with the $\mu_d$ action given by $\zeta \cdot z = \zeta^a z$. It follows that the suspension $\Sigma F_{vw}$ is equivariantly homotopy equivalent to the representation sphere $S^V$ associated to the 1-dimensional complex $\mu_d$-representation $V$ given by $z \mapsto \zeta^a z$. Thus, by equivariant Bott periodicity and the Sebastiani-Thom theorem, we have a weak equivalence
$$K_{\mu_d}(F_f) \xrightarrow{\simeq} K_{\mu_d}(F_{f + vw}),$$
and hence a weak equivalence 
$$K_{\mu_d}(\C^n, F_f) \xrightarrow{\simeq} K_{\mu_d}(\C^{n+2}, F_{f+vw} )$$
(note that, while $F_f$ is not compact, it is equivariantly homotopy equivalent to a compact $\mu_d$-space, so there is no problem with applying equivariant Bott periodicity with $X = F_f$).

Now, set $S':=S[v,w]/(vw)$ and $R':=R[v,w]/(vw)$, where $|v| = a$ and $|w| = b$. Let $\psi: R' \to S'$ be given by setting $u = 0$. Since Kn\"orrer periodicity is induced by tensoring with a fixed matrix factorization, Lemma~\ref{pushtensor} implies that the square
\begin{equation}
\label{kpbp1}
	\xymatrixcolsep{1.1in}
	 \xymatrix{
	 \Ds(S) \ar[d]^-{K} \ar[r]^-{\overline{\varphi^*}} & \Ds(R) \ar[d]^-{K} \\
	\Ds(S')  \ar[r]^-{\overline{\psi^*}}  &  \Ds(R')  }	
\end{equation}
commutes, where the vertical arrows are quasi-equivalences given by Kn\"orrer periodicity. Thus, by iterating Kn\"orrer periodicity sufficiently many times, we may eliminate from Proposition~\ref{hs} the requirement concerning the degree of $f$:

\begin{thm}
\label{main}
There exists a canonical exact triangle of spectra
$$K_{\mu_m}(\C^n, F_f) \to K^{\topp}(\Ds(S)) \xrightarrow{\overline{\varphi^*}} K^{\topp}(\Ds(R))  \xrightarrow{+1}.$$
In particular, there exists a weak equivalence
$$K_{\mu_d}(\C^n, F_f) \xrightarrow{\simeq} K^{\topp}(\Ds(S)).$$
\end{thm}

\begin{rem}
Passing to $K^{\topp}$, Theorem~\ref{main} and (\ref{kpbp1}) yield an isomorphism of exact triangles of spectra of the following form:

\begin{equation}
\label{kpbp}
\xymatrixcolsep{.8in}
		\xymatrix{
		K_{\mu_m}(\C^n,F_f) \ar[r]  \ar[d]^-{\simeq} & K^{\topp}(\Ds(S)) \ar[d]^-{\simeq} \ar[r]^-{\overline{\varphi^*}}  &  K^{\topp}(\Ds(R)) \ar[d]^-{\simeq}    \\
			K_{\mu_m}(\C^{n+2}, F_{f + vw}) \ar[r]  & K^{\topp}(\Ds(S'))  \ar[r]^-{\overline{\varphi^*}} &  K^{\topp}(\Ds(R'))  }
\end{equation}


Observe that (\ref{kpbp}) exhibits a precise sense in which graded Kn\"orrer periodicity and equivariant Bott periodicity are compatible. We refer the reader to Theorem 3.34 of \cite{brown2016knorrer} for a similar compatibility result, at the level of algebraic $K_0$ groups, involving Kn\"orrer periodicity for non-graded matrix factorizations. Note that we have not provided a new \emph{proof} of Bott periodicity, since, as discussed in Section~\ref{ktop}, Bott periodicity is used in the construction of the functor $K^{\topp}$. 
\end{rem}

\subsection{An example: computing a push-forward on topological $K$-theory}
\label{example}
Assume $n>2$ and $w_i$ = 1 for all $i$. Set $q=x_1^2 + \cdots + x_n^2 \in Q$, $S = Q/(q)$, and $R = Q[u] / (q + u)$, where $|u| = 2$. Let $\sZ \subseteq \PP^{n-1}$ denote the projective hypersurface defined by $q$, and let $\sX$ denote the quotient stack $[\Spec(R) \setminus \{0\} / \C^*] \cong \PP^{n-1}$. Let $i: \sZ \into \sX$ denote the closed embedding.

By Theorem~\ref{orlov}, one has semi-orthogonal decompositions
\begin{itemize}
\item[(1)] $\D^b(\coh \sZ) = \langle \OO( -n + 3), \dots, \OO, \Psi(\Ds(S)) \rangle$,
\item[(2)] $\D^b(\coh \sX) =  \langle \OO(-n + 1), \dots, \OO \rangle$.
\end{itemize}
Here, $\Psi$ is the quasi-fully faithful embedding $\Ds(S) \to \D^b(\coh \sZ)$ given by the composition
$$\Ds(S) \xrightarrow{b} \D^b(\gr_{\ge 0} S) \into \D^b(\gr S) \xrightarrow{\widetilde{(-)}} \D^b(\coh \sZ),$$
where $b$ is the functor described in Appendix~\ref{adjoint}. Note that (1) is originally due to Kapranov (\cite{kapranov1988derived} \S4). 

We set
$$K^{\topp}(\sZ):= K^{\topp}(\D^b (\coh \sZ)),$$
$$K^{\topp}(\sX):= K^{\topp}(\D^b (\coh \sX)).$$ The main goal of this section is to compute the map 
$$i_*: K_0^{\topp}(\sZ) \to K_0^{\topp}(\sX).$$
The morphism $i_*: K^{\topp}(\sZ) \to K^{\topp}(\sX)$ of spectra is a special case of a map which appears in diagram~(\ref{eq:relmildiag}) above. We include this example in order to make the abstract computations in the previous sections a bit more concrete.

The most difficult part of this computation is understanding where $i_*$ sends $\Psi(K_0^{\topp}(\Ds(S)))$; this will require some calculations involving graded matrix factorizations of quadrics.

\subsubsection{Push-forward on $K_0^{\alg}$}
\label{pushforward}
The semi-orthogonal decompositions (1) and (2) above, along with Examples~\ref{power} and~\ref{uv} and Kn\"orrer periodicity, imply that one has isomorphisms

\begin{displaymath}
   K^{\topp}_i (\sZ) \cong \left\{
     \begin{array}{ll}
       0 &  i \text{ odd}\\
       \Z^{\oplus (n-1)} &  i \text{ even and } n \text{ odd}\\
       \Z^{\oplus n} & i \text{ even and } n \text{ even}
     \end{array}
   \right.
\end{displaymath}

\begin{displaymath}
   K^{\topp}_i (\sX) \cong \left\{
     \begin{array}{ll}
       0 &  i \text{ odd}\\
       \Z^{\oplus n} &  i \text{ even}\\
     \end{array}
   \right.
\end{displaymath}

Observe also that one has a commutative square

\[
		\xymatrix{
			K_0^{\alg}(\sZ) \ar[d]_{}\ar[r]^-{i_*} & K_0^{\alg}(\sX) \ar[d]_{} \\
			K_0^{\topp}(\sZ)   \ar[r]^-{i_*}  & K_0^{\topp}(\sX) }
\]
where the vertical maps are the forgetful maps. By Examples~\ref{power} and~\ref{uv}, along with Kn\"orrer periodicity, $\Ds(q)$ is quasi-equivalent to either $\Perf (\C)$ or $\Perf (\C) \times \Perf (\C)$, and so the vertical maps are isomorphisms. Hence, it suffices to understand the push-forward on algebraic $K$-theory. 

One has a canonical isomorphism 
$$K_0^{\alg} (\langle \OO(-n+3), \cdots, \OO \rangle) \times K_0^{\alg}(\Ds(S)) \xrightarrow{\text{inc} \times \Psi}  K_0^{\alg}(\sZ),$$
where $\text{inc}$ denotes the map induced by inclusion. The behavior of $i_*$ on the image of $\text{inc}$ is not difficult to understand; for each $l \in \Z$, one has a short exact sequence
$$0 \to \OO_{\sX}(l - 2) \xrightarrow{u} \OO_{\sX}(l) \to i_*(\OO_{\sZ}(l) )\to 0.$$
Thus, $i_*[\OO_{\sZ}(l)] \in K_0^{\alg}(\sX)$ is equal to $[\OO(l )] - [\OO(l-2)]$.

Understanding the image of $\Psi(K_0^{\alg}(\Ds(S)))$ under $i_*$ is a bit more subtle. Let $\C \in \Ds(S)$ denote the complex with the residue field $S/(x_1, \dots, x_n)$ concentrated in degree 0. Using Example~\ref{residuefield}, one sees
$$i_*(\Psi[\C]) = [\OO^{\oplus {2^{n-1}}}] - [\OO(-1)^{\oplus {2^{n-1}}}].$$
To understand where $\Psi \circ i_*$ sends the rest of $K_0^{\alg}(\Ds(S))$, it will be useful to model $\Ds(S)$ by the dg category $\MF^\Z(q)$ of graded matrix factorizations of $q$. By Examples~\ref{power} and~\ref{uv}, we know: 
\begin{itemize}
\item $K_0^{\alg}(\MF^\Z(x_1^2))$ (where $|x_1| = 1$) is free abelian of rank 1, generated by $L:=[\C[x_1] \darrow{x_1}{x_1} \C[x_1](-1)]$.
\item $K_0^{\alg}(\MF^\Z(x_1^2 + x_2^2))$ (where $|x_1| = 1 = |x_2|$) is free abelian of rank 2, generated by $X:=[\C[x_1, x_2](-1) \darrow{x_1 + ix_2}{x_1-ix_2} \C[x_1,x_2](-2)]$ and $X':=[\C[x_1,x_2](-1) \darrow{x_1- ix_2}{x_1 +ix_2} \C[x_1,x_2](-2)]$.
\end{itemize}

Notice that $\Phi(L) = \C$ and $\Phi(X \oplus X') = \C$, where $\Phi$ is the equivalence discussed in Section~\ref{mf}. Thus, letting $K$ denote the Kn\"orrer periodicity functor, we have:

\begin{itemize}
\item if $n$ is odd, the underlying graded free module of $K^{\frac{n-1}{2}}(L)$ has total rank $2^{\frac{n+1}{2}}$, and $\Phi(K^{\frac{n-1}{2}}(L)) \in \Ds(q)$ is a summand of $\C$,

\item if $n$ is even, the underlying graded free modules of $K^{\frac{n-2}{2}}(X)$ and $K^{\frac{n-2}{2}}(X')$ have total rank $2^{\frac{n}{2}}$, and both $\Phi(K^{\frac{n-2}{2}}(X)),\Phi(K^{\frac{n-2}{2}}(X')) \in \Ds(q)$ are summands of $\C$. 
\end{itemize}

It follows that
\begin{itemize}
\item if $n$ is odd, $(i_* \circ \Psi \circ \Phi)(K^{\frac{n-1}{2}}(L)) = [\OO^{\oplus 2^{\frac{n-1}{2}}}] - [\OO(-1)^{\oplus 2^{\frac{n-1}{2}}}] \in K_0^{\alg}(\sX)$,

\item if $n$ is even, $(i_* \circ \Psi \circ \Phi)(K^{\frac{n-2}{2}}(X)) = (i_* \circ \Psi \circ \Phi)(K^{\frac{n-2}{2}}(X') )=  [\OO^{\oplus 2^{\frac{n-2}{2}}}] - [\OO(-1)^{\oplus 2^{\frac{n-2}{2}}}] \in K_0^{\alg}(\sX)$
\end{itemize}

Equip $K_0^{\alg}(\sX)$ with the basis $\{[\OO(-n + 1)], \dots, [\OO] \}$, and equip $K_0^{\alg}(\sZ)$ with the basis
\begin{itemize}
\item $\{[\OO(-n + 3)], \dots, [\OO], [K^{\frac{n-1}{2}}(L)] \}$, if $n$ is odd,
\item $\{[\OO(-n + 3)], \dots, [\OO], [K^{\frac{n-2}{2}}(X)],  [K^{\frac{n-2}{2}}(X')] \}$, if $n$ is even.
\end{itemize}

We have proven:

\begin{prop}
\label{K0}
If $n$ is odd, $i_*: K_0^{\alg}(\sZ) \to K_0^{\alg}(\sX)$ is given by the $n \times (n-1)$ matrix 
$$
\begin{pmatrix}
-1 & 0 & 0 & \cdots & 0 & 0\\
0 & -1 & 0 & \cdots & 0  & 0\\
1 & 0 & -1 & \cdots & 0 & 0 \\
0 & 1 & 0 & \cdots & 0 & 0 \\
0 & 0 & 1 & \cdots & 0 & 0 \\
\vdots & \vdots & \vdots & \vdots & \vdots & \vdots \\
0 & 0 & 0 & \cdots & -1 & 0 \\
0 & 0 & 0 & \cdots & 0 & -2^{\frac{n-1}{2}} \\
0 & 0 & 0 & \cdots & 1 & 2^{\frac{n-1}{2}} \\
\end{pmatrix}
$$

If $n $ is even, $i_*: K_0^{\alg}(\sZ) \to K_0^{\alg}(\sX)$ is given by the $n \times n$ matrix 
$$
\begin{pmatrix}
-1 & 0 & 0 & \cdots & 0 & 0 & 0 \\
0 & -1 & 0 & \cdots & 0 & 0 & 0 \\
1 & 0 & -1 & \cdots & 0 & 0 & 0 \\
0 & 1 & 0 & \cdots & 0 & 0 & 0 \\
0 & 0 & 1 & \cdots & 0 & 0 & 0 \\
\vdots & \vdots & \vdots & \vdots & \vdots & \vdots & \vdots\\
0 & 0 & 0 & \cdots & -1 & 0 & 0 \\
0 & 0 & 0 & \cdots & 0 & -2^{\frac{n-2}{2}} & -2^{\frac{n-2}{2}} \\
0 & 0 & 0 & \cdots & 1 & 2^{\frac{n-2}{2}}& 2^{\frac{n-2}{2}} \\
\end{pmatrix}
$$
\end{prop}

\subsubsection{Topological $K$-theory of the complement}
The complement of $\sZ$ in $\sX$ is $\sU := [ \Spec(\C[x_1, \dots, x_n] / (q - 1)) / \mu_2]$; let $j: \sU \into \sX$ denote the open embedding, and set $K^{\topp}(\sU):= K^{\topp}(\D^b(\coh \sU))$. By the proof of Lemma 3.6 of \cite{halpern-pomerleano:equivariant}, one has an exact triangle
$$K^{\topp}(\sZ) \xrightarrow{i_*} K^{\topp}(\sX) \xrightarrow{j^*} K^{\topp}(\sU) \xrightarrow{+1}$$
of spectra, and hence an exact sequence
$$0 \to K^{\topp}_1(\sU) \to K^{\topp}_0(\sZ) \xrightarrow{i_*} K^{\topp}_0(\sX) \xrightarrow{j^*} K^{\topp}_0(\sU)\to 0.$$ 

One now easily deduces from Proposition~\ref{K0} that
\begin{displaymath}
   K_i^{\topp}(\sU) \cong \left\{
     \begin{array}{ll}
       \Z \oplus (\Z/2^{\frac{n-2}{2}}\Z) &  i \text{ even, } n \text{ even}\\
       \Z &  i \text{ odd, } n \text{ even} \\
       \Z \oplus (\Z/2^{\frac{n-1}{2}}\Z) &  i \text{ even, } n \text{ odd}\\
       0 &  i \text{ odd, } n \text{ odd}\\
      \end{array}
   \right.
\end{displaymath} 

By Theorem 3.9 of \cite{halpern-pomerleano:equivariant}, there exists a canonical weak equivalence $K^{\topp}(\sU) \xrightarrow{\simeq} K_{\mu_2}(F_q)$, where $F_q$ denotes the Milnor fiber of $q$, and $\mu_2$ acts by monodromy. In this case, the action of $\mu_2$ is free (the nontrivial element of $\mu_2$ sends $x \in F_q$ to $-x$), and so we have  
$$K_{\mu_2}(F_q) \cong K(F_q / \mu_2) \cong K(\RR P^{n-1}).$$

Hence, we have recovered the well-known calculation of the topological $K$-theory of real projective space (cf. \cite{atiyah1967k} Proposition 2.7.7).


\section{The Atiyah-Bott-Shapiro construction}
\label{ABS}
We now apply Theorem~\ref{main} to obtain a spectrum-level version of the Atiyah-Bott-Shapiro construction. We begin by recalling the classical Atiyah-Bott-Shapiro construction (\cite{atiyah-bott-shapiro} Part III). For each $n \ge 1$, denote by $C_n$ the complex Clifford algebra associated to the form $q_n:=x_1^2 + \cdots + x_n^2$. Also, set $C_0$ to be the $\Z/2\Z$-graded $\C$-algebra with $\C$ concentrated in degree 0. Let $M(C_n)$ denote the Grothendieck group of the abelian category $\on{mod}_{\Z/2\Z}(C_n)$ of finitely generated $\Z/2\Z$-graded $C_n$-modules. Note that

\begin{displaymath}
   M(C_n) \cong \left\{
     \begin{array}{ll}
        \Z \oplus \Z &  n \text{ even}\\
       \Z  &  n \text{ odd} \\
           \end{array}
   \right.
\end{displaymath} 

The inclusions $i_n: C_n \to C_{n+1}$ induce maps $i_n^*: M(C_{n+1}) \to M(C_n)$
via restriction of scalars. Set $A_n := M(C_n) / i_n^*(M(C_{n+1})).$ The \emph{Atiyah-Bott-Shapiro construction} is the family of canonical isomorphisms
$$A_n \xrightarrow{\cong} \widetilde{K}^0(S^n)$$
for each $n \ge 0$ provided by \cite{atiyah-bott-shapiro} Theorem 11.5. 

Fix $n \ge 1$, and let $R_n$ denote the graded hypersurface ring $\C[x_1, \dots, x_n]/(q_n)$, where $|x_i|=1$ for all $i$. Let $\varphi: R_{n+1} \to R_n$ be the map given by setting $x_{n+1}=0$. The triangle
\begin{equation}
\label{spectrumabs}
K^{\topp}(\Ds(R_n)) \xrightarrow{\overline{\varphi^*}} K^{\topp}(\Ds(R_{n+1})) \to \rK(F_{q_n})
\end{equation}
arising from Theorem~\ref{main} implies that there exists a canonical isomorphism
$$\coker(K_0^{\topp}(\Ds(R_n))  \xrightarrow{\overline{\varphi^*}} K_0^{\topp}(\Ds(R_{n+1}))) \xrightarrow{\cong}  \rK^0(F_{q_n}).$$

Notice that $K_0^{\topp}(\Ds(R_n)) \cong M_{n}$, $K_0^{\topp}(\Ds(R_{n+1})) \cong M_{n-1}$, and $F_{q_n}$ is homotopy equivalent to $S^{n-1}$. The goal of this section is to canonically identify the map $\overline{\varphi^*}: K_0^{\topp}(\Ds(R_n)) \to K_0^{\topp}(\Ds(R_{n+1}))$ with the map $i^*: M_n \to M_{n-1}$, thus demonstrating that the triangle~\ref{spectrumabs} recovers the Atiyah-Bott-Shapiro construction.

We will need the following technical result, which is adapted from Lemma 12.2 of \cite{yoshino1990maximal}:

\begin{lem}
\label{specialform}
Let $Q = \C[x_1, \dots, x_n]$, where $|x_i|$ is a positive integer for each $i$. Let $f \in Q$ be a homogeneous element of degree $2|x_n|$, and suppose $f = g + x_n^2$, where $g \in Q':=\C[x_1, \dots, x_{n-1}] \subseteq Q$. Then there exists an equivalence of additive categories
$$H: [\MF^\Z(f)] \to [\MF^\Z(f)]$$
such that, for each object $F:=(F_0 \darrow{s_0}{s_1} F_1) \in [\MF^\Z(f)]$, 
\begin{itemize}
\item[(1)] there exists a graded free $Q$-module $G_F$ and a degree $|x_n|$ endomorphism $\phi$ of $G_F$ such that 
\begin{itemize}
\item[(a)] $H(F_0 \darrow{s_0}{s_1} F_1) = G_F \darrow{x_n \cdot \id + \phi}{x_n \cdot \id - \phi} G_F(-|x_n|)$,
\item[(b)] letting $e: \gr Q' \to \gr Q$ and $r: \gr Q  \to \gr Q' $ denote the extension/restriction of scalars functors, one has $(e \circ r)(\phi) = \phi$,
\end{itemize}
\item[(2)] there is a natural isomorphism $(F_0 \darrow{s_0}{s_1} F_1) \xrightarrow{\cong} H(F_0 \darrow{s_0}{s_1} F_1)$ in $[\MF^\Z(f)]$.
\end{itemize}
\end{lem}

\begin{proof}
Set $S:=Q/(f)$. There is a well-known equivalence of triangulated categories
$$\coker: [\MF^\Z(f)] \xrightarrow{\cong} \underline{\on{MCM}}^{\gr}(S),$$
where the category on the right is the stable category of graded maximal Cohen-Macaulay $S$-modules. The equivalence is given, on objects, by 
$$(F_0 \darrow{s_0}{s_1} F_1) \mapsto \coker(s_1)$$
and by the evident map on morphisms. 

Set $S':=Q'/(g)$, and let $\overline{r}: \gr S  \to \gr S'$
denote the restriction of scalars functor. Let $F:=(F_0 \darrow{s_0}{s_1} F_1)$ be an object of $[\MF^\Z(f)]$; observe that, by a graded version of the Auslander-Buchsbaum formula, $\overline{r}(\coker(F))$ is a graded free $Q'$-module. Set $G_F:= \overline{r}(\coker(F)) \otimes_{Q'} Q$, and let $\phi: G_F(-|x_n|) \to G_F$ be the map induced, via the functor $\overline{r}(-) \otimes_{Q'} Q$, by $\coker(F)(-|x_n|) \xrightarrow{\cdot x_n} \coker(F)$. We define $H$ to be given, on objects, by 
$$(F_0 \darrow{s_0}{s_1} F_1) \mapsto (G_F \darrow{x_n \cdot \id + \phi}{x_n \cdot \id - \phi} G_F(-|x_n|)),$$ 
and in the obvious way on morphisms. $H$ clearly has properties (1)(a) and (1)(b), and as for (2), it is straightforward to check that the natural map $\coker(F) \to \coker(H(F))$ is an isomorphism. 
\end{proof}

\begin{rem}
The functor $H$ constructed in the proof of Lemma~\ref{specialform} does not preserve shifts. In particular, it is not a triangulated functor.
\end{rem}

\begin{cor} 
\label{restriction}
Let $Q = \C[x_1, \dots, x_n]$, where $|x_i|$ is a positive integer for each $i$. Suppose $f \in Q$ is a homogeneous element of the form $g + x_n^2$, where $g \in \C[x_1, \dots, x_{n-1}] \subseteq Q$. Let $R: \MF^\Z(f) \to \MF^\Z(g)$ be the dg functor given by setting $x_n = 0$. Then one has a commutative triangle
\[
\xymatrixcolsep{.8in}
		\xymatrix{
			K_0^{\alg} \MF^\Z(f) \ar[d]^-{R} \ar[r]^-{T}  &  K_0^{\alg} \MF^\Z(f + u^2)  \\
			K_0^{\alg} \MF^\Z(g)  \ar[ur]^-{K}  }
\]
where $|u| = |x_n|$, $K$ denotes the Kn\"orrer periodicity functor 
$$- \otimes_\MF (\C[x_n,u] \darrow{u + ix_n}{u - ix_n} \C[x_n,u](-|x_n|)): \MF^\Z(g) \to \MF^\Z(f + u^2)$$
(see Remark~\ref{kpint}), and $T$ is the functor
$$- \otimes_{\MF} (\C[u] \darrow{u}{u} \C[u](-|x_n|)): \MF^\Z(f) \to \MF^\Z(f +  u^2).$$ 
\end{cor}

\begin{proof}
Let $F=(F_0 \darrow{s_0}{s_1} F_1)$ be an object in $\MF^\Z(f)$. Choose an equivalence of categories $H: [\MF^\Z(f)] \to [\MF^\Z(f)]$ as in Lemma~\ref{specialform}, and write $H(F) = G \darrow{x_n \cdot \id + \phi}{x_n \cdot \id - \phi} G(-|x_n|)$. Set $G':=G \otimes_Q Q[u]$, and observe that there exists an isomorphism $(K \circ R)(H(F)) \xrightarrow{\cong} T(H(F))$ in $[\MF^\Z(f)]$ given by the pair of maps
$$G' \oplus G' \xrightarrow{\begin{pmatrix} 1 & i \\
			i & 1 \\ \end{pmatrix}} G' \oplus G' \text{, } G'(-|x_n|) \oplus G'(-|x_n|)\xrightarrow{\begin{pmatrix} 1 & i \\
			i & 1 \\ \end{pmatrix}} G'(-|x_n|) \oplus G'(-|x_n|).$$

\end{proof}

Let $[\MF(Q, q_n)]$ denote the homotopy category of non-graded matrix factorizations of $q_n$ (see Definition 2.1 of \cite{dyckerhoff2011compact} for the definition of the dg category of non-graded matrix factorizations). Let $\Delta$ denote the equivalence of categories
$$[\MF(Q, q_n)] \xrightarrow{\cong} \on{mod}_{\Z/2\Z}(C_n)$$
constructed on pages 130-131 of \cite{yoshino1990maximal}; this equivalence was originally established by Buchweitz-Eisenbud-Herzog in \cite{buchweitz1987cohen}.
Let $F: [\MF^\Z(q_n)] \to [\MF(Q, q_n)]$ denote the forgetful functor. Also, as in Section~\ref{pushforward},
\begin{itemize}
\item let $L \in \MF^\Z(x_1^2)$ denote the matrix factorization $\C[x_1] \darrow{x_1}{x_1} \C[x_1](-1)$, and 
\item let $X, X' \in \MF^\Z(x_1^2 + x_2^2)$ denote the matrix factorizations $\C[x_1, x_2](-1) \darrow{x_1 + ix_2}{x_1-ix_2} \C[x_1,x_2](-2)$ and $\C[x_1,x_2](-1) \darrow{x_1- ix_2}{x_1 +ix_2} \C[x_1,x_2](-2)$, respectively.
\end{itemize}
Recall that $K_0^{\alg}\MF^\Z(q_n)$ is a rank one free abelian group generated by $[K^{\frac{n-1}{2}}(L)]$ when $n$ is odd and a rank two free abelian group generated by $[K^{\frac{n-2}{2}}(X)]$ and $[K^{\frac{n-2}{2}}(X')]$ when $n$ is even.

\begin{lem}
\label{BEH}
The map $c: K_0^{\alg}(\MF^\Z(q_n)) \to M_n$ given by
$$[K^{\frac{n-1}{2}}(L)] \mapsto [(\Delta \circ F)(K^{\frac{n-1}{2}}(L))] \text{, if n is odd}$$
$$[K^{\frac{n-2}{2}}(X)] \mapsto [(\Delta \circ F)(K^{\frac{n-2}{2}}(X))] \text{ and } [K^{\frac{n-2}{2}}(X')] \mapsto [(\Delta \circ F)(K^{\frac{n-2}{2}}(X'))] \text{, if n is even}$$
is an isomorphism. Moreover, the diagram
\[
\xymatrixcolsep{.8in}
		\xymatrix{
			K_0^{\alg} (\MF^\Z(q_n))  \ar[d]^-{c}  \ar[r]^-{R}  &  K_0^{\alg} (\MF^\Z(q_{n-1})) \ar[d]^-{c} \\
			M_n  \ar[r]^-{i^*} & M_{n-1} \\
			   }
\]
commutes, where $R$ is as in Corollary~\ref{restriction}.
\end{lem}

\begin{proof}
By a Theorem of Buchweitz-Eisenbud-Herzog, the forgetful functor $F$ is essentially surjective (see Proposition 14.3 of \cite{yoshino1990maximal}); thus, $c$ is surjective. Since $K_0^{\alg}(\MF^\Z(q_n))$ and $M_n$ are free abelian groups of the same rank, $c$ is an isomorphism. The second statement is immediate from the construction of the functor $\Delta$.
\end{proof}

Applying Lemma~\ref{pushtensor}, Corollary~\ref{restriction}, and Lemma~\ref{BEH}, we obtain canonical isomorphisms $M_n \xrightarrow{\cong} K_0^{\topp} \Ds(R_n)$, $M_{n-1} \xrightarrow{\cong} K_0^{\topp} \Ds(R_{n+1})$ making the diagram
\[
\xymatrixcolsep{.8in}
		\xymatrix{
		M_n  \ar[d]^-{\cong} \ar[r]^-{i^*} & M_{n-1} \ar[d]^-{\cong} \\
		K_0^{\topp} \Ds(R_n)   \ar[r]^-{\overline{\varphi^*}}   &  K_0^{\topp} \Ds(R_{n+1}) \\
			   }
\]
commute.

\begin{rem}
\label{absspectra}
The idea of recovering the Atiyah-Bott-Shapiro construction from a localization sequence in $K$-theory is not new; a similar approach is taken by Swan in Section 10 of \cite{swan1985k}. However, our presentation in~\ref{spectrumabs} of the reduced nonconnective topological $K$-theory spectrum of $S^{n-1}$ as the cofiber of a map between two spectra which agrees with the map $i^*: M_n \to M_{n-1}$ upon passing to $\pi_0$ (i.e. a ``lift" of the Atiyah-Bott-Shapiro construction to the level of spectra) appears to be new.
\end{rem}

\appendix
\section{Proof of Theorem~\ref{gradedknorrer}}
\label{kappendix}
\begin{proof}[Proof of Theorem~\ref{gradedknorrer}]
Since $\MF^\Z(f)$ and $\MF^\Z(f + uv)$ are pretriangulated (\cite{ballard2014category} Corollary 3.6), it suffices to show the induced functor
$$K: [\MF^\Z(f)] \to [\MF^\Z(f + uv)]$$
on homotopy categories is an equivalence. Given an additive category $\mathcal{A}$ equipped with an automorphism $T$, let $\mathcal{A}/T$ denote its orbit category, as defined in \cite{keller2005triangulated}. We recall that $\mathcal{A}/T$ has the same objects as $\mathcal{A}$ and morphisms from $X$ to $Y$ given by $\bigoplus_i \Hom_{\mathcal{A}}(X, T^i(Y))$. Recall that $Q$ denotes the graded ring $\k[x_1, \dots, x_n]$, with $|x_i|$ a positive integer for all $i$. Consider the commutative square
$$
\xymatrix{
[\MF^\Z(f)] \ar[r]^-{K} \ar[d]^-{F} & [\MF^\Z(f + uv)] \ar[d]^-{F'} \\
[\MF(Q, f)] \ar[r]^-{\widetilde{K}} & [\MF(Q[u,v], f + uv)],
}
$$
where $[\MF(Q, f)]$ and $[\MF(Q[u,v], f + uv)]$ denote the categories of non-graded matrix factorizations of $f$ and $f + uv$ over the polynomial rings $Q$ and $Q[u,v]$, $F$ and $F'$ are the forgetful functors, and $\widetilde{K}$ is the non-graded Kn\"orrer functor given by tensoring with $\k[u,v] \darrow{u}{v} \k[u,v]$. It is well-known that $\widetilde{K}$ is an equivalence; see, for instance, \cite{orlov2006triangulated}. 

We consider the orbit categories $[\MF^\Z(f)]/(1)$ and $[\MF^\Z(f + uv)]/(1)$, where $(1)$ denotes the automorphism given by twisting. Since $K$ commutes with twists, the above square factors like so:

$$
\xymatrix{
[\MF^\Z(f)] \ar[r]^-{K} \ar[d]^-{\alpha} & [\MF^\Z(f + uv)] \ar[d]^-{\alpha'} \\
[\MF^\Z(f)]/(1) \ar[d]^{\beta} \ar[r]^-{\overline{K}} & [\MF^\Z(f + uv)]/(1)  \ar[d]^-{\beta'} \\
[\MF(Q, f)] \ar[r]^-{\widetilde{K}} & [\MF(Q[u,v], f + uv)],
}
$$
where $\alpha$ and $\alpha'$ are the canonical maps. By Lemma A.7 of \cite{keller2011two}, $\beta$ and $\beta'$ are fully faithful. Hence, $\overline{K}$ is fully faithful. It follows easily that $K$ is fully faithful as well. 

We now show essential surjectivity. Let $X \in \MF^\Z(f + uv)$. Choose $Y \in \MF(Q, f)$ such that $\widetilde{K}(Y) \cong F'(X)$. Let 
$$R: [\MF(Q[u,v], f + uv)] \to [\MF(Q, f)]$$ 
denote the functor which sets $u$ and $v$ to $0$. Clearly $(R \circ \widetilde{K})(Y) = Y \oplus Y[1]$. Thus, 
$$(\widetilde{K} \circ R )(F'(X)) \cong F'(X) \oplus F'(X)[1].$$ 
In particular, $F'(X)$ is a summand of $(\widetilde{K} \circ R )(F'(X))$. Since $\beta'$ is fully
faithful, we may conclude that $\alpha'(X)$ is a summand of $\overline{K}(Z)$ for some $Z \in
[\MF^\Z(f)]/(1)$. We denote the corresponding inclusion by $g: \alpha'(X) \into \overline{K}(Z)$ and
write $g = \oplus_{i \in I} g_i$ for the decomposition of $g$ into its homogeneous components so
that $g_i: \alpha'(X) \to \overline{K}(Z)(i)$ is the image of a morphism in $[\MF^\Z(f + uv)]$. It
follows that $X$ is a summand of $\oplus_{i \in I} K(Z)(i) \cong \oplus_{i \in I} K(Z(i))$ in
$[\MF^\Z(f + uv)]$. Since $K$ is fully faithful, and $[\MF^\Z(f )]$ is idempotent complete
(\cite{ballard2014category} Corollary 3.6), we're done.
\end{proof}

\section{The left adjoint of $\pi \iota : \on{D}^b(\gr_{\ge 0}R) \to \on{D}^{\sg}(R)$}
\label{adjoint}
Let $R = \oplus_{i \ge 0} R_i$ be a graded Gorenstein algebra over a field $\k$. Recall from Section~\ref{orlovsection} the map $\pi \iota : \on{D}^b(\gr_{\ge 0} R) \to \on{D}^{\sg}(R).$ Lemma 2.4 of \cite{orlov2009derived} implies that $\pi \iota$ has a fully faithful left adjoint $b$ that embeds $\on{D}^{\sg}(R)$ as the left orthogonal of $\mathcal{P}_{\ge 0} \into \on{D}^b(\gr_{\ge 0} R)$, where $\mathcal{P}_{\ge j}$ is the smallest thick subcategory of $\on{D}^b(\gr_{\ge 0} R)$ containing $R(e)$ for all $e \le - j$. Burke-Stevenson describe the functor $b$ explicitly in Section 5 of \cite{burke2013derived}; we recall their description here, since we make use of it in Section~\ref{example}.

Given an object $M \in \on{D}^{\sg}(R)$, 
\begin{enumerate}
\item Choose a complex $P$ of graded projective $R$-modules such that
\begin{itemize}
\item $P^i = 0$ for $i \gg 0$,
\item for all $j \in \Z$, there exists $k_j$ such that $P^{-i} \in \mathcal{P}_{\ge j}$ for all $i \ge k_j$, and
\item $P$ is quasi-isomorphic to $M$
\end{itemize}
\item Denote by $P^+$ the subcomplex of $P$ consisting of summands of $R(j)$ with $j \le 0$. Construct a complex $Q$ of graded projective $R$-modules as in (1), this time quasi-isomorphic to $\RHom_{\gr R}(P^+, R)$, where, given graded $R$-modules $L$ and $N$, $$\underline{\Hom}_{\gr R}(L, N) := \bigoplus_{j \in \Z} \Hom_{\gr R}(L(-j), N).$$ 

\item Finally, $b$ sends $M \in \on{D}^{\sg}(R)$ to $\RHom_{\gr R}(Q, R)^-$, the subcomplex of $\RHom_{\gr R}(Q, R)$ consisting of summands of $R(j)$ with $j > 0$. This complex does indeed lie in $\on{D}^b(\gr_{\ge 0} R)$, since $R$ has finite injective dimension as a graded module over itself.
\end{enumerate}

\begin{exm}
\label{residuefield}
Suppose $R$ is a graded hypersurface ring of the form 
$$\k[x_1, \dots, x_n]/(q),$$ where $|x_i|=1$ for all $i$, and $q = x_1^2 + \cdots + x_n^2$. Take $M \in \on{D}^{\sg}(R)$ to be the residue field $R/(x_1, \dots, x_n)$ concentrated in degree $0$. We compute $b(M)$.

\begin{itemize}
\item Step (1): take $P$ to be a minimal graded $R$-free resolution of $M$:
$$  \cdots \to R(-n)^{\oplus m}  \to  R(-n + 1)^{\oplus m} \to R(-n+2)^{\oplus m_{n-2}}\to \cdots \to  R(-1)^{\oplus m_1} \to R \to 0 \to \cdots$$
Notice that the exponents remain constant after $n-1$ steps in the resolution, because $R$ is a hypersurface ring of Krull dimension $n-1$. It is well-known that $m = 2^{n-1}$; see, for instance, \cite{buchweitz1987cohen}.
\item
Step (2): notice that $P = P^+$, in this case. Also, we have
$$\RHom_{\gr R}(M, R) \simeq M(a)[-n + 1],$$
where $a$ is the Gorenstein parameter $n-2$. Thus, we can take $Q$ to be a minimal free resolution of $M(a)$:
$$  \cdots \to R(a - n)^{\oplus m} \to  R(a - (n-1))^{\oplus m} \to \cdots \to R(a-2)^{\oplus m_2} \to R(a-1)^{\oplus m_1} \to R(a) \to 0 \to \cdots $$
where $R(a)$ sits in cohomological degree $-n+1$.

\item Step (3): $\RHom_{\gr R}(Q, R)$ is the complex 
$$ \cdots \to 0 \to R(-a) \to R(-a + 1)^{\oplus m_1} \to \cdots \to R(-a + n-1)^{\oplus m} \to R(-a + n)^{\oplus m} \to \cdots,$$
where $R(-a)$ sits in degree $-n+1$. Hence, $\RHom_{\gr R}(Q, R)^-$ is the complex
$$C:=\cdots \to 0 \to R(1)^{\oplus m} \to R(2)^{\oplus m} \to R(3)^{\oplus m} \to \cdots $$
where $R(1)^{m}$ sits in degree 0. 

Let $\phi: R(-n)^{\oplus m} \to R(-n + 1)^{\oplus m}$ be the map arising in $P$, and set $N:= \coker(\phi)$. The complex with $N(n-2)$ concentrated in degree 0 is, of course, quasi-isomorphic to the twist
$$\cdots \to R(-3)^{\oplus m} \to R(-2)^{\oplus m} \to  R(-1)^{\oplus m}  \to 0 \to \cdots$$
of the tail of $P$, where $R(-1)^{\oplus m}$ sits in degree $0$. Thus, $\RHom_{\gr R}(N(n-2), R) = C$. Since $C$ only has nonzero cohomology in degree 0, we have that 
$$C \simeq \underline{\Hom}_{\gr R}(N(n-2), R)$$ 
(we are using that $N(n-2)$ is maximal Cohen-Macaulay, here). Thus, $b(M)$ is the graded maximal Cohen-Macaulay $R$-module $\underline{\Hom}_{\gr R}(N(n-2), R)$, concentrated in degree 0.
\end{itemize}

One can easily compute the matrix factorization associated to $\underline{\Hom}_{\gr R}(N(n-2), R)$ under the equivalence $\Phi$ of Section~\ref{mf}. In general, given an MCM module $L$, let $L^*$ denote its graded $R$-dual $\underline{\Hom}_{\gr R}(L, R)$. If $L = \coker(F_1 \xrightarrow{s_1} F_0)$, where $F_0 \darrow{s_0}{s_1} F_1$
is a matrix factorization of $q$, then 
$$L^* = \coker((F_0)^*(-2) \xrightarrow{s_1^*} (F_1)^*(-2));$$
that is, $L^*$ corresponds to the matrix factorization
$(F_1)^*(-2) \darrow{s_0^*}{s_1^*} (F_0)^*(-2)$. We apply this to our example. Notice that $N(n-2)$ corresponds to a matrix factorization of the form $Q(-1)^{\oplus 2^{n-1}} \darrow{s_0}{s_1} Q(-2)^{\oplus 2^{n-1}},$
where $Q = \k[x_1, \dots, x_n]$ (recall that $m = 2^{n-1}$). Thus, $\underline{\Hom}_{\gr R}(N(n-2), R)$ corresponds, under $\Phi$, to the matrix factorization
$Q^{\oplus 2^{n-1}} \darrow{s_0^*}{s_1^*} Q(-1)^{\oplus 2^{n-1}}.$
\end{exm}


\bibliographystyle{amsalpha}
\bibliography{Bibliography}

\end{document}